\documentclass[12pt]{article}
\usepackage{amsmath}
\usepackage{amssymb}
\usepackage{amscd}
\usepackage{amsthm}

\usepackage{graphicx}\usepackage{epsfig}
\usepackage{epstopdf}

\newtheorem{dfn}{Definition}[section]
\newtheorem{rem}[dfn]{Remark}
\newtheorem{thm}[dfn]{Theorem}

\newtheorem{lem}[dfn]{Lemma}

\newtheorem{prop}[dfn]{Proposition}

\newtheorem{cor}[dfn]{Corollary}

\newtheorem{example}[dfn]{Example}

\newtheorem{question}[dfn]{Question}

\def\CR{\curvearrowright}
\def\acts{\CR}

\def\R{{\mathbb R}}
\def\eps{\epsilon}
\def\al{\alpha}
\def\be{\beta}
\def\ga{\gamma}

\def\H{\mathbb H}

\def\Ga{\Gamma}
\def\del{\delta}

\def\La{\Lambda}

\def\si{\sigma}

\def\C{{\mathbb C}}

\def\De{\Delta}

\def\embed{\hookrightarrow}

\def\om{\omega}

\def\la{\lambda}
\def\t{\tilde}

\def\<{\langle}
\def\>{\rangle}

\def\about{\asymp}

\begin{document}

\title{RAAGs in Ham}
\author{Michael Kapovich}
\date{\today}


\maketitle

\section{Introduction}

For a graph $\Ga$ let $V(\Ga), E(\Ga)$ denote the vertex and edge sets of $\Ga$.
Let $\Ga$ be a graph with no loops and bigons, i.e., a simplicial complex of dimension $\le 1$. 
Define the Right Angled Artin group (RAAG) $G_\Ga$ with the {\em Artin graph} $\Ga$
by the presentation
$$
\< g_v, v\in V(\Ga)| \quad [g_v, g_w]=1, [vw]\notin E(\Ga)\>.$$
We note that our definition is opposite to the usual one in the theory of RAAGs, 
see e.g. \cite{C}, where one imposes the relators $[g_v, g_w]=1$ for every $[vw]\in E(\Ga)$. However,  our convention 
is in line with the notation in the theory of finite Coxeter groups and Dynkin diagrams. We adopted this notation because it 
is most suitable for the purposes of this paper, while the usual definition leads to heavy notation. 

Given a symplectic manifold $(M,\om)$ we let $Ham(M,\om)$ denote the group of Hamiltonian symplectomorphisms of $(M,\om)$. 
Since, by Moser's theorem, for a closed surface $M$ its symplectic structure is unique up to scaling, we will abbreviate $Ham(M,\om)$ to 
$Ham(M)$ if $M$ is a closed surface. 

Our main result is:

\begin{thm}\label{main}
For every finite $\Ga$ the group $G_\Ga$ embeds in $Ham(S^2)$. Moreover, under this embedding the group $G_\Ga$ 
fixes a closed disk in $S^2$ pointwise.  
\end{thm}

As a corollary of the proof of this theorem we establish the following result proven in the end of the paper: 

\begin{cor}\label{c:main}
For every finite $\Ga$ and every symplectic manifold $(M,\om)$, the group $G_\Ga$ embeds in $Ham(M,\om)$. 
\end{cor}

\begin{cor}
Let $\La\subset O(n,1)$ be an arithmetic lattice of the simplest type, $n\ge 2$. Then a finite index subgroup in $\La$ embeds in $Ham(M,\om)$ for every 
symplectic manifold $(M,\om)$. 
\end{cor}
\proof According to the result of Bergeron, Haglund and Wise \cite{Haglund-Wise}, a finite index subgroup in 
$\La$ embeds in some RAAG $G_\Ga$. Now, the result follows from Corollary \ref{c:main}. \qed 

\medskip 
In contrast, suppose that $M$ is a closed oriented surface of genus $\ge 1$ with area form $\om$. 
Then it was proven first by  L.~Polterovich \cite{Polterovich02} and, later, by Franks and Handel \cite{FH} using different methods, 
that every irreducible nonuniform arithmetic group $\La$ of rank $\ge 2$ does not embed in $Ham(M,\om)$. 
Furthermore, Franks and Handel \cite{FH} extended this result to certain nonuniform rank 1 lattices, e.g., lattices in $PU(2,1)$.


\medskip 
{\bf Outline of the proof.} Theorem \ref{main} is proven in three steps.

{\bf Step 1}. Let $M$ be a closed connected oriented surface to which $\Ga$ embeds. 
For technical reasons, it will be convenient to assume that $M$ is not the torus. We first prove

\begin{thm}\label{step1}
The group $G_\Ga$ embeds in $Ham(M)$.  
Moreover, each Artin generator $g_v$ of $G_\Ga$ acts on $M$ as an ``iterated Double Dehn twist''  
$\Psi(g_v)$ supported in a homotopically trivial annulus in $M$.
\end{thm}

The key to verifying injectivity of $\Psi: G_\Ga\to Ham(M)$ is that the action $G_\Ga\acts M$ preserves a certain finite subset 
$P\subset M$, so that the restriction $G_\Ga\acts M'=M\setminus P$ projects to a faithful  representation to the mapping class group of $M'$,  
$G_\Ga\to Map(M')$. Faithfulness of this representation follows from a special case of a theorem of L.~Funar \cite{Funar} 
(see also the more recent papers by T.~Koberda \cite{Koberda} and by M.~Clay, C.~Leininger and J.~Mangahas \cite{CLM}). This part of our paper is similar to the arguments by J.~Crisp and B.~Wiest \cite{Crisp-Wiest2}. 

\medskip 
{\bf Step 2} (Lifting).  If $\Ga$ were planar,  Theorem \ref{step1} would imply Theorem \ref{main}. 
In general, of course, $\Ga$ need not be planar (or even admit a finite planar orbi-cover, see \S \ref{orbicovers}), 
however, it has a planar universal cover (e.g., the disjoint union of simplicial trees). 
Suppose, therefore, that $M$ has genus $\ge 2$. Then we lift the action $\Psi: G_\Ga\acts M$ to the universal cover $\tilde{M}$ of $M$, 
which we identify with the hyperbolic plane, i.e., the unit disk $D$ in $S^2=\C\cup \{\infty\}$. We let $\om_0$ be the Euclidean area form on an open disk containing $D$; extend $\om_0$ smoothly to an area form $\om_0$ on $S^2$.

Let $D'=D-P'$ denote the punctured disk where $P'$  is the preimage of $P$ in $D$. 
Let $Ham(M,P)$ denote the subgroup of $Ham(M)$ fixing $P$ pointwise. We have an (injective) homomorphism
$$
\iota: Ham(M,P)\to Ham(D,P')
$$
obtained by choosing an appropriate lifting of Hamiltonian diffeomorphisms.  
We thus obtain the lift $\tilde\Psi=\iota\circ \Psi$ of the homomorphism $\Psi$. Then we show that $\tilde\Psi$ projects injectively to 
the mapping class group $Map(D')$.

Each generator $g_v$ of $G_\Ga$ acts (via $\t\Psi$) on $D$ as a product of infinitely many commuting $N$-iterated 
Double Dehn twists preserving the hyperbolic area form. However, $\t\Psi(G_\Ga)$, of course, does not preserve $\om_{0}$. 
Then we modify each of the Double Dehn twists in the product decomposition of $\t\Psi(g_v)$ to obtain a new diffeomorphism $\rho_0(g_v)$ 
which is isotopic to $\t\Psi(g_v)$ on the punctured disk $D'$ and is the time-$N$-map for the appropriately chosen function $H_v: D\to \R$ with respect to $\om_{0}$. It then follows that the resulting representation 
$$\rho_0: G_\Ga\to Ham(D, \om_{0})$$
is again faithful. We will see that for each $v$, $H_v$ extends by zero to a $C^{1,1}$-function on $S^2$ and 
$\rho_0(g_v)$ extends Lipschitz-continuously (by the identity) to the entire sphere, so we can think of it as a Lipschitz Hamiltonian symplectomorphism. However, the function $H_v$ need not be $C^2$-smooth and $\rho_0(g_v)$ need not even be differentiable.

\medskip 
{\bf Step 3} (Approximation). 
The last step of the proof is an approximation argument: We approximate $H_v: S^2\to \R$ by a mollifier,  a smooth function 
$\eta_{\eps} H_{v}$ which depends analytically on $\eps>0$ and converges to $H_v$ uniformly on compacts in the open disk $D$ as $\eps\to 0$. Each function $\eta_{\eps} H_{v}$ determines its own time-$N$ map $\rho_\eps(g_v)$ and we obtain an analytic family of representations $\rho_\eps: G_\Ga\to Ham(S^2)$, $\eps>0$, which converge to $\rho_0$ as $\eps\to 0$. Then we establish that the representations $\rho_\eps$ are injective for all but countably many $\eps>0$, thereby proving Theorem \ref{main}.


\medskip
{\bf Questions.}

\begin{question}
Is it true that every RAAG $G_\Ga$ admits a quasi-isometric embedding $G_\Ga\to Ham(S^2)$ with the $L_2$-metric on $Ham(S^2)$? See 
\cite{Crisp-Wiest2} for the case of RAAGs with planar $\Ga$. 
\end{question}

\begin{question}
Let $M$ be a closed oriented surface and  let $\La$ be a K\"ahler  group which is not virtually a surface group and is not virtually abelian. Does $\La$ embed in $Ham(M)$? Does there exist an infinite group with property T that 
embeds in $Ham(M)$? Conjecturally, the latter question has negative answer, see \cite{Fisher}. 
\end{question}

\begin{question}
(A.~Berenstein) Is it true that every Artin group (not necessarily right-angled) embeds in $Ham(S^2)$?  
\end{question}

\begin{question}
(V.~Kharlamov) Which RAAGs embed in $Diff(S^1)$? 
\end{question}

\medskip
{\bf Acknowledgments.} The author was supported by the NSF grant DMS-09-05802, he is also grateful to Max Plank Institute for Mathematics in Bonn for its hospitality. The author is grateful to the referee for useful remarks, to Thomas Koberda from whom he learned about the results of \cite{Koberda}, to 
Pierre Py for pointing out at \cite{Fisher} and to Mark Sapir and Leonid Polterovich for numerous discussions. This work was motivated by the problem which arose during the Oberwolfach Workshop ``Geometric Group Theory, Hyperbolic Dynamics and Symplectic Geometry" in 2006 and the author thanks MFO for hosting the workshop.

\section{Orbi-covers of graphs and embeddings of RAAGs}\label{orbicovers}

This section is purely algebraic and its goal is two-fold: The key result is Lemma \ref{diag} which will allow us to 
apply the results of L.~Funar for proving Theorem \ref{step1}. The point is that Funar's theorem \cite{Funar} 
(and, similarly, results of T.~Koberda \cite{Koberda}) deals with subgroups of Mapping Class groups generated by 
(iterated) Dehn twists, while we will be using (iterated) Double Dehn twists. We then prove Lemma \ref{cover} which is a vast 
generalization of Lemma \ref{diag}.  
Lemma \ref{cover} allows, in the case of graphs which admit planar emulators, to avoid the analytical arguments in Parts 2 and 3 
of the proof of Theorem \ref{main}. This could be useful if 
one were to construct explicit embeddings of various RAAGs  in $Ham(S^2)$ (the proof of Theorem \ref{main} is non-constructive). 
Lemma \ref{cover} could be also useful for solving the classification problem of embeddings between RAAGs (compare \cite{Kim-Koberda}).

\medskip
The {\em double} $D\Ga$ of a graph $\Ga$ is defined as follows. Start with the disjoint union $\Ga\times \{-1, 1\}$ of two copies of $\Ga$. 
Then for every edge $[v,w]$ of $\Ga$ we add edges $[v\times i, w\times j]$, $i\ne j$, to  $\Ga\times \{-1, 1\}$. The result is $D\Ga$. We 
will use the notation $v^\pm:=v\times \pm 1$ for $v\in V(\Ga)$.

If $G_\Ga$ is the RAAG with the Artin graph $\Ga$, we call $G_{D\Ga}$ the {\em double} of $G_\Ga$. 
Then we have the diagonal homomorphism
$$
\del: G_\Ga \to G_{D\Ga}, \quad \del(g_v)=g_{v^+} g_{v^-}.
$$
Note that since $[v^+, v^-]\notin E(D\Ga)$, the order in the product $g_{v^+} g_{v^-}$ is irrelevant.  
To see why the map of the generators of $G_\Ga$ extends to a group homomorphism, suppose that $[v,w]\notin E(\Ga)$, i.e., the 
generators $g_v, g_w$ commute. Then, by the definition of $D\Ga$, the vertices $v^\pm, w^\pm$ are not connected by edges in $D\Ga$. 
Hence, $\del(g_v)$ and $\del(g_w)$ commute as well.

\begin{lem}\label{diag}
The homomorphism $\del: G_\Ga \to G_{D\Ga}$ is injective.
\end{lem}
\proof We have the natural projection $\pi:  G_{D\Ga} \to G_{\Ga}$, $\pi(g_{v^+})=g_v, \pi(g_{v^-})=1$. Then, clearly,
$\pi\circ \del=id$. Hence, $\del$ is injective. \qed

This lemma is a special case of a more general result on embeddings of RAAGs proven below. 

An {\em orbi-cover} (or a {\em branched-cover}, also known as a {\em weak cover}) 
of a graph is a map of graphs $p: \De\to \Ga$ which  is a locally-surjective
graph morphism. (Local surjectivity means that for every vertex $v\in \De$ and every edge $e$ of $\Ga$ incident to $p(v)$, 
there exists an edge $j$ of $\De$ which is mapped to $e$ by $p$. In other words, the $p$ is locally surjective as a map of 
topological spaces. A map of graphs is a graph morphism if it sends  vertices to vertices and edges homeomorphically to edges.)  
A {\em planar emulator} of a graph $\Ga$ is a finite orbi-cover $p: \De\to \Ga$ with planar $\De$, see e.g. 
\cite{Hlineny}. For instance, the graphs $K_5, K_6$ are not planar but admit planar 2-fold finite covers, since they embed in $\R P^2$. 
(Recall that $K_n$ is the complete graph on $n$ vertices.)  Moreover, there are finite graphs $\Ga$ which admit finite planar emulators 
but admit no finite planar covers \cite{Rieck}. 

I am grateful to Yo'av Rieck for the following example:

\begin{example}
Suppose that $\Ga$ is a finite 1-dimensional simplicial complex 
where every vertex $\Ga$ has valence $\ge 6$. Then $\Ga$ does not admit planar emulators. Indeed, let us show first that $\Ga$ is not planar itself. An embedding  $\Ga\to S^2$ would define a 
cell complex decomposition of $S^2$. Triangulating each 2-cell from a vertex, 
we obtain a new cell decomposition where each 2-cell is a triangle and every vertex has valence $\ge 6$. Let $v, e, f$ denote the number of vertices, edges and faces of this decomposition, where $2e=3f$, $2e\ge 6v$. Then the Euler characteristic 
computation yields:
$$
2=\chi(S^2)= v -e + f = v - \frac{1}{3} e \le 0 . 
$$
Contradiction. Suppose that $p: \De\to \Ga$ is a finite orbi-cover of $\Ga$ as above. Let $\Delta'$ denote the subgraph of $\Delta$, which is the maximal simplicial complex in $\Delta$ containing all the vertices of 
$\Delta$. Every vertex of $\Delta'$ still has valence $\ge 6$. Thus, $\Delta'$ and, hence, $\Delta$, cannot be planar. 
\end{example}
 
Given a RAAG $G_\Ga$ and a finite orbi-cover  $p: \De\to \Ga$, one defines a ``diagonal'' homomorphism $\del=p^*: G=G_\Ga\to \t{G}=G_\De$ by
$$
p^*(g_v):= \prod_{x\in p^{-1}(v)} g_x.
$$
Note that all the generators $g_x, x\in p^{-1}(v)$ of the group $G_\De$ commute (since $x, y\in p^{-1}(v)$ are never connected by an edge in $\De$). It is immediate that $p^*: G_\Ga\to G_\De$ is indeed a homomorphism: If $v, w\in V(\Ga)$ and $[vw]\notin \Ga$ then for any $x\in p^{-1}(v), y\in p^{-1}(w)$, 
$[xy]\notin \Ga$. 

\begin{lem}\label{cover}
$\del: G\to \t{G}$ is injective. 
\end{lem}
\proof We will use the normal forms for the elements of RAAGs. We first order the vertices of $\Ga$; we lift this order to a lexicographic order on $V(\De)$. 
Then a normal form of $g\in G$ is the product of generators 
$$
w=g_{v_1}^{\pm 1}... g_{v_k}^{\pm 1}
$$
with the condition that the word $w$ contains no subwords of the form $g_v g_v^{-1}$ and if $[uv]\in E(\Ga)$, $u<v$ and $g_v^{\pm 1}$ precedes  $g_u^{\pm 1}$ in $w$, then between these letters  in $w$ there is a letter $g_z^{\pm 1}$ such that $[zv]\in E(\Ga)$. Then every $g\in G$ admits a normal form and this normal form is unique. The reduction process of a word $w$ to the normal form is as follows:

A pair of consecutive letters $g_{y}^{\pm 1} g_x^{\pm 1}$ in $w$ is an {\em inversion} if $y>x$ and $[xy]\notin E(\Ga)$. Then, in order to reduce $w$ to its normal form use the commutation relation to reduce the number of inversions (``shuffling'') and cancel appearances of the products $g_v g_v^{-1}$ (``cancellation''). We refer the reader to \cite{Hermiller-Meier} for the details. 

Suppose now that $w=g_{v_1}^{\pm 1}... g_{v_k}^{\pm 1}$ is a normal form for $g\in G$. The image $\del(w)\in \t{G}$ defined as 
$$
\del(g_{v_1}^{\pm 1})... \del(g_{v_k}^{\pm 1})
$$
need not be in normal form. We claim, however, that the length of the normal form of $\del(w)$ is the same as the length of $\del(w)$, i.e., 
no cancellations in the reduction process occur. 
Indeed, the only way we can get a cancellation is that $w$ (or $w^{-1}$) contains
$$
w=... g_v ... g_v^{- 1}...
$$
Then, since $w$ is a normal form, between these appearances of $g_v$ and $g_{v}^{-1}$ there is some $g_u$ (or $g_u^{-1}$) 
so that $[uv]\in E(\Ga)$:
$$
w=... g_v ... g_u ... g_v^{-1}...
$$ 
(We assume that the sub-word between $g_v $ and $g_v^{-1}$ is the shortest where a cancellation in $\del(w)$ is possible.) 
Lifting $w$ to $\t{G}$, we see that for each $x\in p^{-1}(v)$ and $y\in p^{-1}(u)$ there exists an edge $[xy]\in E(\De)$. Therefore, 
shuffling the generators of $\t{G}$ would not allow us to move any $g_y^{-1}$ appearing in the lift of $g_v^{-1}$ past the lift of 
$g_u$. Therefore, we would be unable to cancel any of these $g_y^{-1}$ (in the lift of $g_v$) with any $g_y$ (in the lift  of $g_v$). Thus 
$\del$ is injective and, moreover, is a quasi-isometric embedding $G\to \t{G}$. \qed 

\begin{rem}
Let $D\Ga$ be the double of $\Ga$. Then the natural map $D\Ga\to \Ga$ is a 2-fold orbi-cover and Lemma \ref{diag} is a corollary of 
Lemma \ref{cover}. 
\end{rem}

\medskip 
We now observe that, if $\Ga$ admits a planar emulator $p: \De\to \Ga$,
then Theorem \ref{step1} implies that $G_\Ga$ embeds in $Ham(S^2)$, even though, $\Ga$ need not be planar. Indeed, for such $\Ga$ and the 
planar emulator $\Delta\to \Ga$ we would have an embedding $G_\Ga\embed G_\Delta$ (by Lemma \ref{cover}). Without loss of generality we may assume that $\De$ is a simplicial complex. Theorem \ref{step1} then shows that $G_\Delta$ embeds in $Ham(S^2)$. By composing the two embeddings, we obtain an embedding  $G_\Ga\embed Ham(S^2)$. 
Thus, for $G_\Ga$ so that $\Ga$  admits a planar emulator, Steps 2 and 3 of the proof of Theorem \ref{main} are not needed.

\section{Basic facts of surface topology and hyperbolic geometry}\label{prelim}

{\bf Mapping class group.} Let $M$ be a connected oriented surface, possibly with boundary ($M$ need not be compact). 
Let $Homeo(M, \partial M)$ denote the group of heomeomorphisms of $M$ fixing the boundary pointwise and let 
$Homeo_0(M, \partial M)$ be the identity component of $Homeo(M, \partial M)$.  
Then the {\em Mapping Class group} of $M$ is the group
$$
Map(M):= Homeo(M, \partial M)/ Homeo_0(M, \partial M). 
$$  
We will use the notation $[f]$ for the projection of $f\in Homeo(M, \partial M)$ to $Map(M)$.

Recall that the  group of outer automorphisms, $Out(\Pi)$, of a group $\Pi$ is the quotient $Aut(\Pi)/Inn(\Pi)$, where $Inn(\Pi)$ consists of inner automorphisms of $\Pi$. We will use the notation $[\phi]$ for the projection of $\phi\in Aut(\Pi)$ to $Out(\Pi)$.  
Given a surface $M$ as above, we have a 
natural homomorphism $\nu: Map(M)\to Out(\Pi)$, $\Pi=\pi_1(M)$ defined as follows. 
If $f\in Homeo(M)$ had a fixed point $x\in M$, then $\nu([f])$ would be defined 
as the projection of the induced map $f_*: \pi_1(M,x)\to \pi_1(M, x)$. 
In general, one uses instead the induced map $f_*: \pi_1(M,x)\to \pi_1(M, f(x))$, 
where $x\in M$ is a base-point. Choosing a path $\zeta$ in $M$ connecting $x$ to $f(x)$ and attaching the appropriate ``tail'' 
to the loops based at $f(x)$, one obtains a map $f_{\bullet}: \pi_1(M,x)\to \pi_1(M, x)$. The choice of $\zeta$ is, of course, not 
canonical, so  $f_{\bullet}$ is not well-defined. However, projecting to $Out(\Pi)$ eliminates the ambiguity and one, thus, obtains the homomorphism $\nu: Map(M)\to Out(\Pi)$, see e.g. \cite[\S 2.9]{Ivanov}. This homomorphism, in general, is neither surjective 
not injective. However, by a theorem usually attributed to Baer, Dehn and Nielsen, if $M$ is a closed surface, 
then  $\nu$ is an isomorphism. Moreover, if $M$ has empty boundary then $\nu$ is injective. 
See e.g.  \cite[\S 2.9]{Ivanov} and references therein or \cite[\S 8]{Farb-Margalit}. 

\medskip 
{\bf Hyperbolic plane.} 
In what follows, we will be using the Poincar\'e model of the hyperbolic plane $\H^2$, i.e., the unit disk $D\subset \C$ with the 
metric
$$
ds^2= \frac{4|dz|^2}{(1-|z|^2)^2}. 
$$
We will also regard $D$ as a disk in the 2-sphere $S^2$ which is the $1$-point compactification of the complex plane. 
The boundary circle $S^1$ of $D$ is the {\em circle at infinity} of $\H^2$. In this model, the group of orientation-preserving isometries 
$Isom_+(\H^2)$ of $\H^2$ is the group of linear-fractional transformations stabilizing $D$, thus, $Isom_+(\H^2)\subset PSL(2,\C)$. 
The subgroup $Isom_+(\H^2)$ of $PSL(2,\C)$ consists of linear-fractional transformations of the form:
$$
\si(z)=e^{i\theta}\frac{z-a}{-\bar{a}z+1},  \quad |a|<1, a=-e^{-\theta}\si(0). 
$$

\medskip
{\bf Mapping class group and homeomorphisms of $S^1$.} 
Suppose now that $M'$ is a surface without boundary which admits a complete hyperbolic metric of 
finite area which we fix from now on. Set $\Pi':=\pi_1(M')$. Lift the hyperbolic metric on $M'$ to the universal cover of $M'$. 
Then the latter is a complete simply-connected surface of curvature $-1$; therefore, it is isometric to the hyperbolic plane $\H^2$. 
Using this isometry we identify the universal cover with $\H^2=D$, the hyperbolic plane. With this identification, the group $\Pi'$ 
is identified with the group of covering transformations of the universal cover $\H^2\to M'$. Then, $\Pi'$ becomes a discrete 
subgroup of $Isom_+(\H^2)\subset PSL(2,\C)$. Since the surface $M'$ has finite area, the {\em limit set} of $\Pi'$ is the entire circle $S^1$, see e.g. \cite[Theorem 12.15]{Rat}. 

\medskip 
Let $f: M'\to M'$ be a homeomorphism of $M'$, $\nu([f])=[\rho]\in Out(\Pi')$, where  
$\rho\in Aut(\Pi')$. As above, we identify $\Pi'$ with a subgroup of $Isom_+(\H^2)$. 
 Let $\tilde{f}$ denote a lift of $f$ to the hyperbolic plane, the universal cover of $M'$. 
(The lift is unique up to postcomposition with 
elements $\ga\in \Pi'$.) Then one of the lifts $\tilde{f}$ is $\rho$-equivariant (different choices of lifts yield maps which are equivariant under automorphisms $Inn(\ga)\circ \rho$, where $Inn(\ga)$ are the inner automorphisms of $\Pi'$ induced by some $\ga\in \Pi'$). In particular, $\tilde{f}$ admits a homeomorphic extension $h=h_\rho: S^1\to S^1$, where $S^1$ is the boundary circle of the hyperbolic plane (the Poincar\'e disk). Moreover, $h$ is again $\rho$-equivariant and depends only on $\rho$ (and not on $f$). 
Furthermore, $h=Id$ iff $\rho=Id$. 
We refer the reader to \cite{Casson} or \cite[\S 8]{Farb-Margalit} for proofs of these results. 

Thus, we obtain a map 
\begin{equation}\label{map}
\rho\mapsto h_\rho. 
\end{equation}
It is elementary to verify that this map projects to a  map of quotients 
$$
[\rho]\mapsto [h_\rho], \quad Out(\Pi')\to  \Pi' \backslash Homeo(S^1), 
$$
where $[h_\rho]$ denotes the coset $\Pi' h_\rho$.

\medskip
{\bf Restrictions of automorphisms of surface groups to normal subgroups.} 
As an application of the correspondence (\ref{map}) we obtain: 

\begin{lem}\label{normal}
Let $\La'\triangleleft \Pi'$ be a nontrivial normal subgroup of $\Pi'$. Let $\rho\in Aut(\Pi')$ be an element induced by 
some homeomorphism of $M'$, so that $\rho(\La')=\La'$. Then $[\rho]\in Out(\Pi')\setminus \{ 1\}$ implies that  
$\rho|\La'$ projects to a nontrivial element of $Out(\La')$. 
\end{lem}
\proof Suppose that $[\rho|\La']=1\in Out(\La')$. Without loss of generality, we can assume that $\rho|\La'=Id$ (otherwise, we 
replace $\rho$ with a suitable composition $Inn(\la)\circ \rho$, where $\la\in \La'$).  
As before, we realize $\Pi'$ and its subgroup $\La'$ as discrete subgroups of $Isom(\H^2)$ acting on $D$, and, 
hence, on its boundary circle $S^1$. Since $\La'$ is normal in $\Pi'$ and $\La'\ne 1$, it follows that the limit set of $\La'$ is 
the same as the limit set of $\Pi'$, i.e., is the entire circle $S^1$ (see \cite[Theorem 12.1.16]{Rat}). 
The automorphism $\rho\in Aut(\Pi')$ is induced by a homeomorphism $h=h_\rho: S^1\to S^1$. 
Since $\rho$ fixes all elements of $\La'$, it follows from the equivariance condition 
$$
\ga\circ h= \rho(\ga)\circ h= h \circ  \ga
$$
that $h$ fixes all fixed points of all nontrivial elements $\ga\in \La'$. These fixed points are dense in the limit set of $\La'$ 
(see \cite[Theorem 12.1.7]{Rat}). Therefore, $h$ fixes the limit set of $\La'$ pointwise. Since the limit set of $\La'$ is the entire $S^1$, it follows that $h=Id$. Hence, $\rho$ is a trivial automorphism of $\Pi'$ as well. \qed 

\medskip
We will use the following special case of Lemma \ref{normal}. Let $M$ be a closed connected oriented surface of genus $\ge 2$. 
Let $P\subset M$ be a nonempty finite subset and set $M':=M\setminus P$. Pick a point $x\in  M'$. Set $\Pi:=\pi_1(M), \Pi':=\pi_1(M')$. 
We then equip $M$ with a hyperbolic metric. 
Let $p: \t{M}\to M$ be the universal cover. As before, this allows us to identify $\t{M}$ with $D=\H^2$ and the fundamental group 
$\Pi$ with a subgroup of $Isom_+(\H^2)\subset PSL(2,\C)$. We set $D':=p^{-1}(M')\subset D$, then $D'$ is a disk with infinitely many punctures at the points of $P':=p^{-1}(P)$. Let $\t{x}$ be a lift of $x$ to $D'$. 
The restriction $p':=p|D' : D'\to M'$ is again a regular covering map 
(with the group $\Pi$ of deck-transformations). Therefore, by the basic covering theory,  
$\La:=\pi_1(D',\t{x})$ projects isomorphically to a (nontrivial) normal subgroup in $\Pi'=\pi_1(M')$ 
with the quotient group $\Pi= \Pi'/ p'_*(\La)$. Let $f: M\to M$ be a homeomorphism preserving $P$ and fixing the point 
$x\in M'$. Then $f$ lifts uniquely to a homeomorphism 
$\tilde{f}: D\to D$ preserving $D'$ and fixing $\t{x}$.

\begin{cor}\label{cor:normal}
If $f|M'$ is not isotopic to the identify, then   $\tilde{f}: D'\to D'$ is not isotopic to the identity either. 
\end{cor}
\proof Let $[\rho]\in Out(\Pi')$ be induced by $f$, i.e., $\nu([f])=[\rho]$. 
In particular, $\rho$ preserves the subgroup $\La':=p'_*(\La)$ and the map $\tilde{f}: D'\to D'$ induces $[\rho|\La']$ in the sense that 
the following diagram is commutative:
$$
\begin{array}{ccc}
\! \La~~ &\stackrel{\tilde{f}_*}{\longrightarrow} &  \! \La~~ \\
\downarrow p'_* & ~ & \downarrow   p'_*  \\
\! \Pi' ~~&\stackrel{\rho={f}_*}{\longrightarrow} & 	\!   \Pi' ~~\\
\end{array}
$$  
By Lemma \ref{normal}, $[\rho|\La']\in Out(\La')$ is nontrivial. Therefore,  $\tilde{f}_*\in Out(\La)\setminus \{1\}$. 
Hence,  $\tilde{f}$ is not isotopic to the identity. \qed 


\section{Proof of Theorem \ref{step1}}

{\bf Hamiltonian symplectomorphisms.} Let $(M,\om)$ be a symplectic manifold, $H: M\times \R\to \R$ a smooth function. Using the form $\om$ one then converts the differential form $dH(x,t)$ to a time-dependent vector field $X_H(x,t)$ on $M$:
$$
\om(X_H, \xi)=dH(\xi), \quad \xi\in TM.
$$
Consider the ODE
\begin{equation}\label{eq:ham}
\frac{\partial F(x,t)}{\partial t} = X_H(x,t) 
\end{equation}
on the manifold $M$. Solutions $f_t(x):=F(x,t)$ of this equation are {\em Hamiltonian symplectomorphisms} of the manifold $(M,\om)$. 
Note that for a non-compact manifold $M$ the ODE (\ref{eq:ham}) may not have solutions defined on the entire $M$ for any $t>0$. 
In general, Hamiltonian symplectomorphisms of $(M,\om)$ form a pseudo-group $Ham(M,\om)$. However, if $M$ is closed,  
$Ham(M,\om)$ is a group. 

\medskip 
{\bf Double Dehn Twists.} Recall that smooth  manifolds with boundary satisfy {\em Moser's Lemma}:

\begin{thm}
\cite[Lemmata 1 and 2]{Greene-Shiohama}
Suppose that $M$ is a smooth manifold with boundary and 
$\om, \om'$ are volume forms so that
$$
\int_M \om = \int_M \om'.
$$
Then there exists a diffeomorphism $f$ isotopic to the identity, which carries $\om$ to $\om'$. Moreover, if $\om'$ varies continuously in $C^\infty$ topology, $f$ can be also chosen to  vary continuously in $C^\infty$ topology. 
\end{thm}

In particular, if $(M, \om)$ is a symplectic surface and $A\subset M$ is an annulus with smooth boundary, then $(A,\om)$ is
symplectomorphic to a product annulus $A_a:=S^1\times [-a, a]$ with the product area form, where $S^1$ is the unit circle.

We define a {\em twist} Hamiltonian symplectomorphism
$f: A_a\to A_a$ as follows.  Pick a smooth function $H(s,t)=h(t)$, so that $h$ vanishes
(with all its  derivatives) at $-a$ and $a$, and $h'(0)=2\pi$. Let $X_H$ be the associated Hamiltonian vector field on
$A_a$. The field $X_H$ is constant with respect to the $s$-coordinate and tangent to the circles $A^1\times t, t\in [-a,a]$.
Let
$f: A_a\to A_a$ be the corresponding  time-$1$ Hamiltonian symplectomorphism:
$$
f(x)=F(x, 1), \quad F(x,0)=id, \quad F=F(x, \tau), \quad \frac{\partial F}{\partial \tau} = X_H.
$$  
We let $A^\pm_{a}:= S^1\times [\pm a, 0]$ denote the subannuli in $A_a$ with the common boundary circle $C=S^1\times 0$, which we will 
call the {\em central circle} of $A_a$. We let $f_\pm$ denote the restrictions $f|A^\pm_{a}$
extended by the identity to the rest of $A_a$. Lifting $H, X_H$ and $f$ to the universal cover $\tilde{A}_a$
of $A_a$ we see that the lift of $f$ fixes the boundary lines of $\tilde{A}_a$ and acts on the line
$\R\times 0$ (the lift of the central circle $C$) as the translation 
by $2\pi$. Therefore, both $f_\pm$ are Dehn twists on $A_a$ and $f_-$ is isotopic to the inverse of
$f_+$ relative to the boundary of $A_a$. By abusing the terminology, we will say that $f$ is the {\em rotation by $2\pi$ along $C$.} 

\medskip 
Let $(M, \om)$ be a symplectic surface, $A\subset M$ be a smooth annulus which is symplectomorphic to some $A_a$.
We will use the notation $C_A$ (the {\em central circle} of $A$) for the circle in $A$ corresponding to $C=S^1\times 0\subset A_a$.   Using the
symplectomorphism $A_a\to A$ we carry the maps $f, f_\pm: A_a\to A_a$, function $H$ and the Hamiltonian vector field $X_H$ to maps
$f, f_\pm: A\to A$, function $H$ and vector field $X_H$ on $A$. The maps $f, f_\pm: A\to A$
extend by the identity to the rest of the surface $M$. We will use  the notation $f_A, f_{A,\pm}$ for the extensions.  
Then $f_A: M\to M$ is a smooth symplectomorphism and $f_A=f_{A,+}\circ f_{A,-}$. The maps $f_{A,\pm}$ are Dehn twists on $M$ which are, up to
isotopy, inverses to each other. Moreover, $f_A$ is a Hamiltonian symplectomorphism since the above function $H: A\to \R$ and its Hamiltonian
vector field $X_H$ extend by zero to the rest of $M$. 

Pick a point $p\in C_A$, then the map $f_A: M\setminus p \to  M\setminus p$ has infinite order in the
mapping class group of this punctured surface, provided that the annulus $A$ is {\em essential},
i.e., each component of $M \setminus A$ has negative Euler characteristic. We will refer to the map
$f_A: M\to M$ as a {\em Double Dehn twist} (such maps are also known as {\em point-pushing maps}).

\medskip
{\bf Construction of homomorphisms of RAAGs to $Ham(M, \om)$.}

Let $G_\Ga$ be a RAAG with the Artin graph $\Ga$. Since $\Ga$ is finite, there exists a 
closed oriented surface $M$ which admits an embedding $j:\Ga\to M$. Without loss of generality, we may assume that $M$ is not the torus. We equip $M$ with an area form $\om$. There exists a collection of closed disks
${\mathcal B}:=\{B(v): v\in V(\Ga)\}$ so that:

1. Each $B(v)$ has smooth boundary.

2.  $B(v)\cap B(w)\ne \emptyset$ (for $v\ne w$) iff $[vw]\in \Ga$.

3. Whenever $B(v)\cap B(w)\ne \emptyset$, their boundary circles $C_v, C_w$ intersect transversally and in exactly two points.  

4. Triple intersections of discs are empty.

\begin{rem}
One can construct such ${\mathcal B}$ as follows: Let $\Ga'$ be the barycentric subdivision of $\Ga$.  
For each $v\in V(\Ga)$ take a sufficiently small smooth disk neighborhood $B(v)\subset M$ of $j(Star(v))$, where
$Star(v)$ is the star of $v$ in $\Ga'$.
\end{rem}

\begin{figure}[htbp]  
   \centering
   \includegraphics[width=3in]{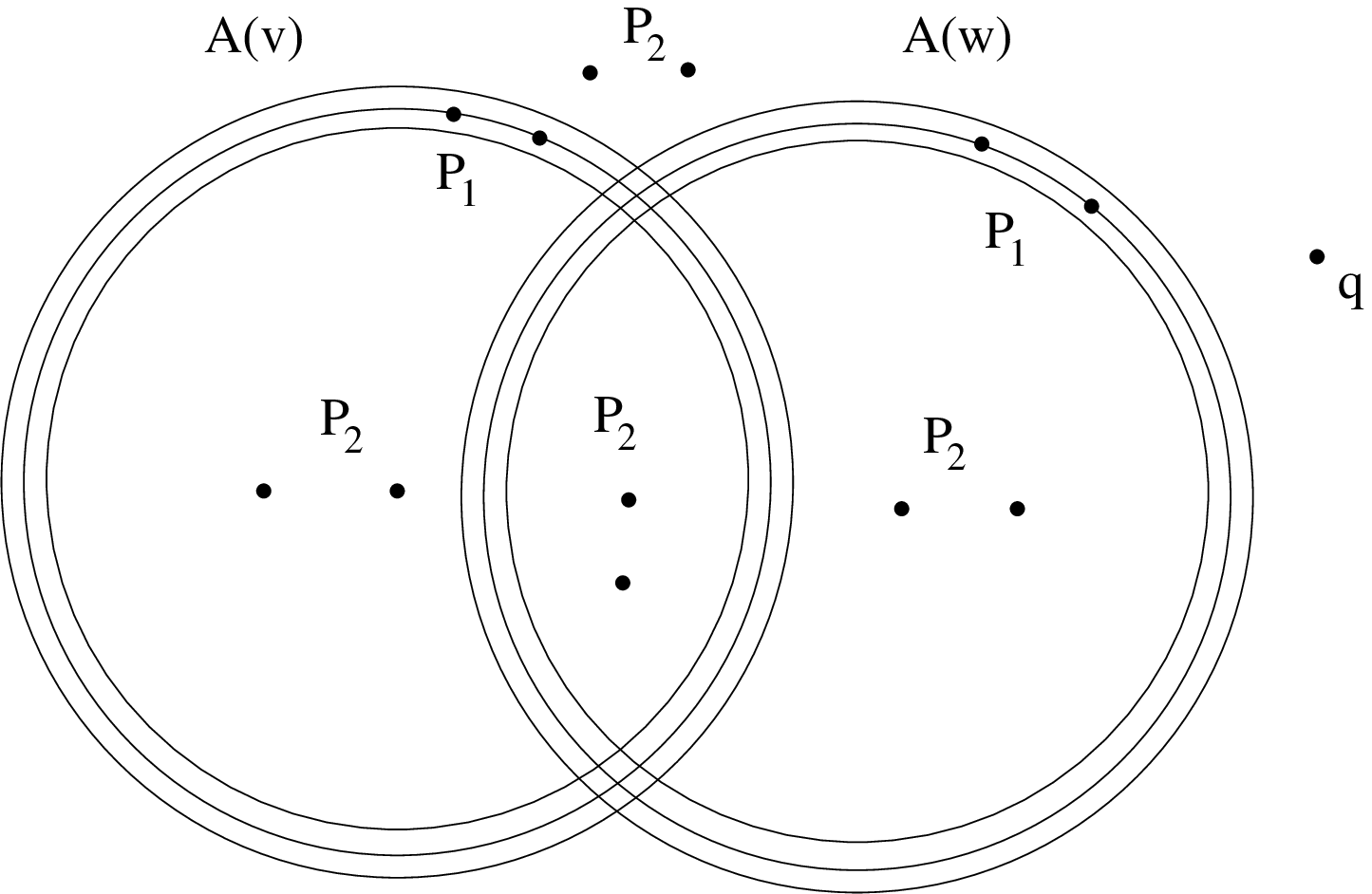} 
   \caption{}
   \label{F0}
\end{figure}

\medskip
Then we thicken each circle $C_v$ to an annulus $A(v)$ in such a way that the nerve of the resulting collection of annuli $\{A(v):  v\in V(\Ga)\}$ is still isomorphic to $\Ga$ and the annuli intersect as in Figure \ref{F0}. We identify each annulus $A(v)$ with the corresponding symplectomorphic product annulus $A_a$, where $a$ depends on $v$. Accordingly, we carry all the notation introduced for $A_a$ to the annulus $A(v)$.  We will identify the circles $C_v$ in the above construction with the central circles $C_{A(v)}$ of the annuli $A(v)$. We note that all annuli $A(v)$ are inessential in $M$, since $M\setminus A(v)$ contains the disk
$D(v)\setminus A(v)$.

We now define a certain finite subset $P\subset M$ fixed by all the double Dehn twists in the annuli $A(v)$. The points of $P$ will serve as punctures on $M$. 
For each vertex $v\in V(\Ga)$ we pick a 2-element set $P_v\subset C_v$,  contained in the connected component of 
$$
C_v \setminus \bigcup_{w\in V(\Ga), w\ne v} A(w). 
$$   
Set
$$
P_1:=\bigcup_{v\in V(\Ga)} P_v 
$$
We let $P_2$ denote a subset of 
$$
M\setminus \bigcup_{v\in V(\Ga)} A(v). 
$$
containing two points in each component of this surface. Lastly, 
pick some $q\in M\setminus (P_1\cup P_2)$ which does not belong to any of the annuli 
$A(v)$ and any of the disks $B(w)$.  We set 
$$
P:= P_1\cup P_2 \cup \{q\} 
$$
and let $M':=M\setminus P$. See Figure \ref{F0}. Now, all the circles in 
$$
\bigcup_{v\in V(\Ga)} \partial A(v)$$
are essential and pairwise non-isotopic in $M'$.

We let $f_v:=f_{A(v)}, f_{v,\pm}:=f_{A(v),\pm}$ denote the Hamiltonian Double Dehn twists and Dehn twists of $(M,\om)$
determined by the symplectic annuli $(A(v), \om)$ and the functions $H_{A(v)}: A(v)\to \R$ corresponding to the function $H: A\to \R$.  
It is clear from the construction that each Dehn twist $f_{v,\pm}, v\in V(\Ga)$, fixes the set $P$ pointwise.
Moreover, each $f_v|M'$ is isotopic to the product of two commuting (isotopically) nontrivial Dehn twists  $f_{v,\pm}: M'\to M'$. 
Let $\al_v^\pm$ denote the boundary circles of $A(v)$ corresponding to the circles $S^1\times \{\pm a\}\subset A_a$. 
Recall that the Dehn twists in an annulus $A\subset M$ is determined, up to isotopy, by the isotopy class of one of the boundary circles 
$\al$ of $A$. The corresponding  element of the Mapping Class group is called the {\em Dehn twist along $\al$}. Therefore, we will think 
of the isotopy classes of maps $f_{v,\pm}$ as Dehn twists along $\al_v^\pm$. 

Since the boundary circles $\al_v^\pm$ of the annuli $A(v)$ are essential and non-isotopic in $M'$,  
it follows that the maps $f_{A(v),\pm}$ are pairwise non-isotopic on $M'$ and, moreover, generate distinct cyclic subgroups of 
the mapping class group $Map(M')$. 

Since $A(v)\cap A(w)=\emptyset$ for $[vw]\notin E(\Ga)$, $f_v$ commutes with $f_w$ whenever $[vw]\notin E(\Ga)$.
It follows that the map $\psi: g_v\to f_v$ determines a homomorphism
$$
\psi: G_\Ga \to Ham(M,\om).
$$
The image of $\psi$ is contained in $Diff(M, P)$, the subgroup of $Diff(M)$ fixing $P$ pointwise.
Moreover, for each natural number $N$ we have a homomorphism 
$$
\psi_N: G_\Ga \to Ham(M,\om), \quad \psi_N(g_v)=f_v^N.  
$$
The diffeomorphism $f_v^N$ is the time-$N$ map of the Hamiltonian $H_{A(v)}$. 

\medskip 
{\bf Funar's Theorem.} In what follows we will need a theorem of L.~Funar \cite[Theorem 1.1]{Funar} formulated below. 

Let $S$ be a compact oriented surface with at least one boundary component. We will use the notation  $F$ 
for a noncompact surface obtained from $S$ by attaching a punctured disk with at least two punctures to each boundary circle of $S$. 
(The number of punctures will be specified later on.) We observe that $Map(S)$ injects in $Map(F)$ 
(see e.g. \cite[Theorem 2.7.I]{Ivanov}).

Let ${\mathcal A}:=\{a_1,...,a_m\}$ be a system of simple closed oriented loops on $S$. 
We require that these loops have the least intersection number in their isotopy classes. 
One says that the system of loops ${\mathcal A}$ is {\em sparse} if for some choice of paths $\ga_i$ connecting $q$ to $a_i$, 
the loops 
$$
b_i:= \ga_i^{-1} a_i \ga_i 
$$
based at $q$ generate a free subgroup of rank $m$ in $\pi_1(S,q)$. We next note that a simple sufficient 
condition for a system of loops ${\mathcal A}$ to be sparse is that they define a linearly independent system of elements of $H_1(S)$.  
(or, equivalently, of $H_1(F)$). 
Here and below, we use homology with real coefficients. Indeed, since $\pi_1(S)$  is free, the group 
generated by the loops $b_i$ is necessarily free. Its rank equals the rank of the subspace in $H_1(S)$ spanned by 
the elements $[a_i]$.

We now assume that ${\mathcal A}$ is {\em sparse} in $S$. Let $D_{a_i}$ denote the Dehn twist (right or left) in $a_i$. Define the RAAG 
$G_\La$, where $\La$ is the incidence graph of the collection of loops ${\mathcal A}$, i.e., $V(\La)=A$, $[a_i, a_j]\in E(\La)$ iff $a_i\cap a_j\ne \emptyset$. Let $q$ be a 
point in the interior of $S$, disjoint from the curves in ${\mathcal A}$. 

\begin{thm}\label{thm:funar}
[L.~Funar] Under the above conditions, for every $N\ge 2$, the natural homomorphism  $\phi_N: G_\La\to Map(S\setminus \{q\})$, 
$\phi_N: g_{a_i}\mapsto [D^N_{a_i}]$, is injective. In particular, the homomorphism  $G_\La\to Map(F)$ obtained by composing 
$\phi_N$ with the embedding $Map(S)\embed Map(F)$ is injective as well. 
\end{thm} 

 We will apply this theorem in the case of punctured surfaces as follows. 

\begin{prop}
For $N\ge 2$ the homomorphism $\Psi:=\psi_N: G_\Ga \to Ham(M,\om)$ is injective.
\end{prop}
\proof Clearly, it suffices to show that the composition
$$
G_\Ga \stackrel{\psi_N}{\to} Diff(M') \stackrel{\pi}{\to} Map(M')
$$
is injective.  
We let $D\Ga$ denote the double of $\Ga$ and $G_{D\Ga}$ be the corresponding double Artin group. 
Then we have natural homomorphisms
$$
\phi_N: G_{D\Ga}\to Diff(M'), \quad \phi_N(g_{v^\pm})= f_{v^\pm}^N.
$$
Then $\psi_N=\phi_N\circ \delta: G_\Ga\to Map(M')$, where $\delta: G_\Ga\embed G_{D\Ga}$ is the
diagonal embedding as in Lemma \ref{diag}. We let $\bar\phi_N: G_{D\Ga}\to Map(M')$ and 
$\bar\psi_N: G_{\Ga}\to Map(M')$ denote the compositions $\pi\circ \phi_N$ and $\bar\phi_N\circ \delta$. 
We observe that the homomorphism $\bar\phi: G_{D\Ga}\to Map(M')$ has the property that 
each Artin generator $g_{v^\pm}$ of $G_{D\Ga}$ maps to the Dehn twist along the boundary curve $\al_v^\pm$ of $A(v)$. 
We claim that for every $N\ge 2$ the homomorphism $\bar\phi_N: G_{D\Ga}\to Map(M')$ is injective. In view of Lemma \ref{diag}, 
this would imply injectivity of $\bar\psi_N$, and, hence, of $\psi_N$ as well. 

\medskip 
We will derive injectivity of $\bar\phi_N$ from Funar's theorem above. 
We define a compact surface $S$, as in Funar's theorem, as follows. 
Set
$$
A:= \bigcup_{v\in V(\Ga)} A(v)
$$
and   define the compact surface $T$   
$$
T:=M \setminus \bigcup_{v\in V(\Ga)} int(B(v)). 
$$
Recall that for every $v\in V(\Ga)$ the set $P_1\cap C_v$ is a 2-element subset $P_v$ contained in a connected component of 
$(M\setminus A)\cup A(v)$. 
We let $\be_v\subset C_v$ be the arc connecting the points of $P_v$ which is disjoint from $A\setminus A(v)$. Clearly, 
all the  Dehn twists $f_{A(w),\pm}$ 
fix every arc $\{\be_{v}\}$. 
 
Recall that the set of punctures $P_2$ in $M$ contains some points in 
$T$. For each component $T_j$ of $T$ we pick a disk $U_j\subset int(T_j)$ containing $T_j\cap P_2$ and not containing the point $q$. 
Thus, each $U_j$ contains at least two points of $P_2$. Set 
$$R:=A\cup T \setminus \bigcup_{j} int(U_j).$$ 
Lastly, cut $R$ open along the arcs $\be_{v}$ defined above and let $S$ denote the resulting surface. See Figure \ref{F3}.

\begin{figure}[htbp] 
   \centering
   \includegraphics[width=3in]{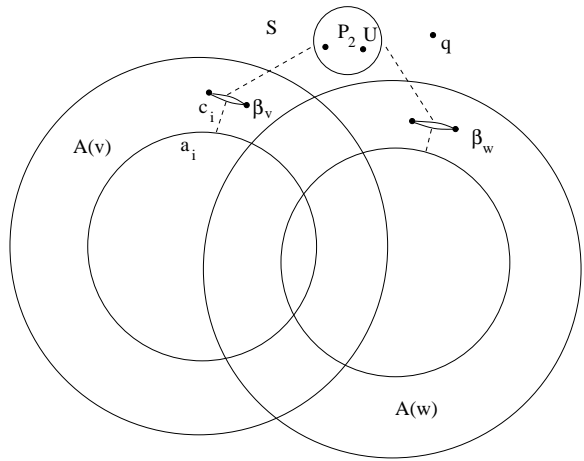} 
   \caption{In this example, $S$ is the sphere with $6$ holes.}
   \label{F3.eps}
\end{figure}

We orient the loops $\al_{v^\pm}$ in an arbitrary fashion. 
We claim that the system ${\mathcal A}$ of curves $\al_{v^\pm}$ in $S$ is linearly independent in $H_1(S)$. 
Indeed, for each loop $a_i:=\al_{v^\pm}$ there exists a properly 
embedded arc $c_i\subset S$ which intersects $a_i$ in exactly one point (possibly, the end-point of $c_i$),  
and intersects the boundary circle corresponding to $\be_v$ at its end-point. 
Moreover, the arc $c_i$ is disjoint from all curves in ${\mathcal A}$ different from $a_i$, see Figure \ref{F3}. 
The relative cycles $[c_i]\in H_1(S, \partial S)\cong H^1(S, \partial S)$ 
are then Poincar\'e dual to $[a_i]\in H_1(S)$, which implies linear independence of    
${\mathcal A}$  in $H_1(S)$. 

\medskip 
Lastly, we observe that the surface $F:=M'$ is obtained from $S$ by attaching punctured disks with at least two punctures each and that the point $q\in P$ belongs to $S$. Therefore, by Theorem \ref{thm:funar}, the homomorphism 
$$
\bar\phi_{N}: G_{D\Ga} \to Map(S\setminus \{q\}) \embed Map(M')
$$ 
is injective.  We conclude that $\psi_N$ is injective as well. \qed 

\begin{rem}
As an alternative to the above argument, one could use the results of \cite{Koberda}, which, however, do not provide an explicit estimate  
on $N$.  
\end{rem}

This finishes the proof of Theorem \ref{step1}. \qed

\section{Lifting to the universal cover}

We continue with the notation introduced in the previous section. Thus, we have a 
closed connected oriented surface $M$ with $\chi(M)<0$ and a punctured 
surface $M'=M\setminus P$, where $\emptyset \ne P\subset M$, and $P$ is finite. 
We have a collection of Hamiltonian diffeomorphisms $f_v\in Ham(M)$ supported on annuli $A(v)\subset M$. Each $f_v$ fixes the set $P\cup \{q\}$ and also fixes the {\em central circle} $C_v\subset A(v)$ pointwise. The incidence graph of the collection of annuli $A(v)$ is the graph $\Ga$ and we have homomorphisms 
$$
\bar{\psi}_N: G_\Ga \stackrel{\psi_N}{\to} Ham(M')\subset Homeo(M,P) \stackrel{\pi}{\to} Map(M').   
$$    
The homomorphism ${\psi}_N$ send Artin generators $g_v$ of $G_\Ga$ to iterated Double Dehn twists $f_v^N$ and 
the homomorphisms $\bar{\psi}_N$ are injective for all $N\ge 2$. We also pick a point $x\in M'$ fixed by all the maps $f_v$.

\medskip 
Our next goal is to lift the Double Dehn twists $f_v$ to the universal cover of the surface $M$. Without loss of 
generality, we may assume that $M$ has genus $\ge 2$, i.e., it admits a hyperbolic structure which we fix from now on. 
We let $\om$ denote the area form of the hyperbolic metric and $\t\om$ its lift to $\t{M}$. 
As in \S \ref{prelim} we identify the universal cover $\t{M}$ of $M$ with the hyperbolic 
plane $\H^2$ embedded in $\C$ as the unit disk $D$. In particular, the area form $\t\om$ is the hyperbolic area form 
$$
\t\om= \frac{dx \wedge dy}{(1-(x^2+y^2))^2}. 
$$

The universal cover $p: D\to M$ yields a covering map $D'\to M'$, where 
$D'=p^{-1}(D)$. Let $\t{x}$ denote a lift of $x$ to $D'$. Then every $f\in \psi_N(G)$ admits a unique lift $\tilde{f}$ fixing $\t{x}$. 
We then obtain homomorphisms $\tilde{\psi}_N: G_\Ga\to Diff(D)$ which send every $g\in G_\Ga$ to the diffeomorphism 
$\tilde{f}: D\to D$. 
Since $\bar{\psi}_N$ is injective (for $N\ge 2$), by applying 
Corollary \ref{cor:normal} we conclude that for every $g\in G_\Ga \setminus \{1\}$, the homeomorphism 
$$\tilde{\psi}_N(g): D'\to D'$$ 
is not isotopic to the identity. By the construction, each $\tilde{\psi}_N(g)$ preserves the area 
form $\t\om$ on $D$. Furthermore, since each annulus $A(v)$ is null--homotopic in $M$,  its preimage  
$p^{-1}(A(v))$ is a disjoint union of annuli in $D$. Each map $\tilde{\psi}_N(g_v)$ is the product of 
commuting $N$-th iterates of Double Dehn twists supported  in the annular components of $p^{-1}(A(v))$. To describe these maps 
more explicitly, let $\tilde{H}_{A(v)}$ denote the lift of the function $H_{A(v)}$ to $D$. Then each 
$\tilde\psi_1(g_v)$ is the time-$1$ Hamiltonian map with respect to the function   $\tilde{H}_{A(v)}$  and the symplectic structure $\t\om$. 
Hence, each $\t\psi_N$ sends $G_\Ga$ to $Ham(D, \t\om)$

In view of the above observations, each homomorphism $\t\psi_N$ (for $N\ge 2$) is injective and, moreover, projects to an injective homomorphism 
$$
G_\Ga\to Map(D'). 
$$
The problem, however, is that the symplectic structure $\t\om$ does not extend to a symplectic structure on the entire sphere $S^2$. 
We, therefore, have to replace it with a symplectic structure $\om_{0}$ on $S^2$, which restricts to the Euclidean area form on an open neighborhood of the closure of $D$. Our next goal is to replace $\t\psi_N(g_v)$ with another infinite product of commuting iterated Double Dehn twists which are Hamiltonian with respect to $\om_0$.  Of course, this will also require correcting the functions $\tilde{H}_{A(v)}$ on $D$. 

\medskip 
{\bf Correcting the functions $\tilde{H}_{A(v)}$.}
For each annulus $A(v)\subset M$ we choose its (homeomorphic) lift to $D$ which we will again denote ${A}(v)\subset D$. 
For each $\si\in \Pi$ we let $\om_{v,\si}$ be the symplectic form on $A(v)$ defined by taking the pull-back of $\om_0$ via
$$
\si: A(v)\to \si(A(v))
$$
and then rescaling by some $\la_{v,\si}^{-2}\in \R_+$, so that
$$
\int_{A(v)} \om_{v,\si}=1.
$$
Clearly,
\begin{equation}\label{zero}
\lim_{\ell(\si)\to\infty} \la_{v,\si}^2=0,
\end{equation}
where $\ell$ is a word metric on $\Pi$. The constants $\la_{v,\si}$ are $\asymp |\si'(z)|, z\in A(v)$, where  
$$
a\asymp b \iff a=O(b) \hbox{~~and~~} b=O(a). 
$$

\begin{lem}
The forms $\om_{v,\si}$ form a precompact set  in $C^{\infty}$ topology.
\end{lem}
\proof Since $\si\in PSL(2,\C)$, we have
$$
\om_{v,\si}= |\theta_{v,\si}(z)| \om_0
$$
where $\theta_{\si}(z):=\theta_{v,\si}(z)=\la_{v,\si}^{-2} \si'(z)^2$ is a function holomorphic
on ${A}(v)$ and having unit $L^1$-norm. Such functions $C^\infty$ subconverge to a holomorphic function on 
the closed annulus $A(v)$. Indeed, each $|\theta_\si|$  is harmonic and $L_1$-norms of these functions are uniformly bounded on $A(v)$. 
Thus, by the mean value property for harmonic functions, $C^0$-norms of the functions $|\theta_\si|$ are again uniformly bounded. Therefore, the holomorphic functions $\theta_{\si}$ form a normal family and, hence, by Cauchy integral formula, $C^\infty$--subconverge to a holomorphic function.  

\medskip 
One can make the above argument more explicit as follows. Without loss of generality, we may assume that the interior of 
$A(v)$ contains $0\in \C$ (otherwise we replace $\Pi\subset Isom_+(\H^2)$ with its conjugate via an element of $Isom_+(\H^2)$ sending 
$0$ to an interior point of $A(v)$). Each linear-fractional transformation $\si=\si(z)$ has the form
$$
\si(z)=e^{it}\frac{z-a}{-\bar{a} z+1},  \quad t=t_\si\in [0, 2\pi], a=a_\si=-e^{-t_\si}\si(0), |a|<1.
$$
Hence,
$$
\si'(z)= \frac{(1-|a_\si|^2)e^{it_\si} }{(-\bar{a}_\si z+1)^2}. 
$$
Since $\la_{v,\si} \asymp |\si'(0)|$, 
$$
\la_{v,\si} = Const_{\si} |\si'(0)|,  
$$
where $Const_{\si}>0$ is bounded away from $0$ and $\infty$. After passing to a subsequence, we obtain: 
$$
\lim a_\si = b, \hbox{~~where~~} |b|=1, \quad \lim Const_{\si}= Const, \quad \lim t_\si =s\in [0,2\pi],  
$$
Here and below all limits are taken with respect to the word norm $\ell(\si)$   diverging to infinity. 
Therefore, 
$$
\lim \theta_\si(z)= \lim Const_{\si} \frac{\si'(z)^2}{(1-|a_\si|^2)^2}=$$ 
$$
\lim Const_{\si} \frac{e^{it_\si} }{(-\bar{a}_\si z+1)^2}= 
Const \frac{e^{2is} }{(-\bar{b} z+1)^4}. 
$$
Moreover, the convergence is uniform on $A(v)$ since $A(v)$ is compact in $D$. \qed 

\medskip 
We retain the notation ${C}_v$ for the circle in ${A}(v)\subset D$ which covers the central circle $C_v\subset A(v)\subset M$.
The circle ${C}_v$  divides the hyperbolic area of the annulus ${A}(v)$ in half but this
need not be the case with respect to the form $\om_{v,\si}$. Nevertheless, by the above compactness 
lemma in conjunction with Moser's lemma, we can choose a $C^\infty$--precompact family of area-preserving diffeomorphisms 
$$
({A}(v), \om_{v,\si})\to A=A_{1/2}=S^1\times [-1/2, 1/2]$$
which carry the circle $C_v$ to the 
round circles $S^1\times b_{v,\si}\subset A$, where $b_{v,\si}$ form a precompact subset of the open annulus $S^1\times (-1/2, 1/2)$. 
We now repeat the construction of Hamiltonian Double Dehn twists on the annulus, except we will insist on having a rotation by $2\pi$ along the circles
$S^1\times b_{v,\si}$ instead of $S^1\times 0$. To this end, we will be using hamiltonians $\hat{H}_{v,\si}: A\to \R$
so that $\hat{H}_{v,\si}(s,t)=h(t)$, $h'(b_{v,\si})=2\pi$.
Pull-back these functions to the annuli $({A}(v), \om_{v,\si})$. We obtain a $C^\infty$--precompact family of functions on $A(v)$.  
The corresponding Double Dehn twists on ${A}(v)$ will rotate ${C}_v$ by $2\pi$.  

Note however that $\si^*(\om_0)=\la_{v,\si}^{2}\om_{v,\si}$ and, hence, we cannot use the above hamiltonians to define Double Dehn twists with respect to
the forms $\si^*(\om_0)$ since the resulting time-1 maps would rotate ${C}_v$ by $2\pi\la_{v,\si}^{-2}$. Therefore, the correct family of functions
$\t{H}_{v,\si}: A(v)\to \R$ is given by the pull-back of
$$
\la_{v,\si}^{2} \hat{H}_{v,\si}
$$
via the symplectomorphisms $({A}(v), \om_{v,\si})\to A$. Clearly, the functions $\t{H}_{v,\si}$ converge to zero in $C^\infty$ topology
on the annulus $A(v)$.

We now define the function ${H}_{v,\si}: \si(A(v))\to \R$ by $\t{H}_{v,\si}\circ \si^{-1}$. Every such function defines (with respect to 
the form $\om_0$) a Hamiltonian Double Dehn twist 
$f_{\si(A(v))}$ on the annulus  $\si(A(v))$ which is isotopic (rel. boundary and the punctures on $C_v$)
to the the Double Dehn twist $\tilde\psi(g_v)$ restricted to $\si(A(v))$. 

We define $H_v: S^2\to \R$ by 
\begin{equation}\label{H_v}
H_v|_{\si(A(v))}= H_{v,\si}, \quad H_v(z)=0 \hbox{~~for} \quad z\in S^2\setminus \bigcup_{\si\in \Pi} \si(A(v)). 
\end{equation}
Accordingly, we extend the maps $f_{\si(A(v))}$ by the identity on the complement of $\si(A(v))$ in $S^2$ and use the notation $f_v$ for the product 
of the resulting commuting double Dehn twists:
$$
f_v=\prod_{\si\in \Pi} f_{\si(A(v))}. 
$$ 
Since 
$$
\lim_{\ell(\si)\to \infty} diam(\si(A(v)))=0,
$$
it is clear that $f_v: S^2\to S^2$ is a homeomorphism.  Since
$$
\lim_{\ell(\si)\to \infty} \|H_{v,\si}\|_{C^0}=0,
$$
it follows that $H_v: S^2\to \R$ is continuous. Clearly, $f_v|D$ is smooth and is the time-1 map of $H_v|D$ with respect to $\om_0$. 
Moreover, $[f_u, f_v]=1$ provided that $[uv]\in E(\Ga)$, since the support sets of the maps $f_u, f_v: D\to D$ are disjoint.  

Choosing $N\ge 2$, we thus obtain a homomorphism
\begin{equation}\label{rho_0}
\rho_0: G_\Ga \to Ham(D,\om_0), \quad \rho_0(g_v)=f_v^N. 
\end{equation}
Since each $f_v$ is isotopic to $\t\psi_N(g_v)$ on $D'$, it follows that the homomorphisms 
$$
\t\psi_N, \rho_0: G_\Ga\to Diff(D')$$ 
have the same projection to $Map(D')$. Since the projection of $\t\psi_N$ to $Map(D')$ was 1-1, it follows that $\rho_0$ also projects injectively. In particular, $\rho_0$ is 1-1 as well. 

\medskip 
Our next goal is to analyze smoothness of the functions $H_v$ and maps $\rho_0(g_v)$. 
The following lemma (and its corollary) is not needed for the proof of Theorem \ref{main} and we include the proof only for the sake of completeness and as a warm-up for the proof of Lemma \ref{fast} which will play an important role in smoothing the functions $H_v$.

\begin{lem}\label{lip}
For every $v\in V(\Ga)$, $H=H_v: S^2\to \R$ is $C^{1,1}$-smooth, i.e., it has Lipschitz differential.
\end{lem}
\proof We only have to verify smoothness on the boundary circle $S^1$ of the unit disk $D$. The function $H$ on the annulus
$A_{v,\si}$ equals $\la_{v,\si}^{2} \hat{H}_{v,\si}\circ \si^{-1}$. Since the derivative of $\si$ on $A(v)$ is of the order of $\la_{v,\si}$, we conclude
that  $dH|_{\si(A(v))}$ converges uniformly to zero as $\ell(\si)\to\infty$. Moreover, the second derivatives of $H$ are uniformly bounded
(by the upper bound on the $C^2$-norm of $\hat{H}_{v,\si}$).
It remains to check that $dH$ vanishes at the boundary points $\xi\in S^1$ of the unit disk. Observe that the Euclidean 
distance from the annulus $A_{v,\si}$ to $S^1$ 
is  $\about \la_{v,\si}$. Therefore, if $d(\xi, z)= R$, where $z\in A_{v,\si}$, then $R\ge C_1\la_{v,\si}$; this implies that
$$
\frac{H(z)}{R}\le C_1 C_2\la_{v,\si} 
$$
where 
$$
\|\hat{H}_{v,\si}\|_{C^0}\le C_2.$$
Thus (\ref{zero}) implies that $H$ has vanishing derivative on $S^1$. The statement that $dH$ is Lipschitz on the closed disk follows
from the above bound on the 2-nd derivative. \qed

\begin{cor}
For each $g\in G_\Ga$, $\rho_0(g)$ is Lipschitz on $S^2$.  
\end{cor}

\begin{cor}
The homomorphism $\rho_0: G_\Ga\to Homeo(S^2)$ is injective. Its image consists of bilipschitz symplectomorphisms of 
$(S^2,\om_0)$ which are Hamiltonian with respect to $C^{1,1}$ functions on $S^2$.    
\end{cor}

This proves a version of Theorem \ref{main} but with very low regularity of symplectomorphisms of $S^2$. Our goal is to replace 
these   bilipschitz symplectomorphisms with infinitely differentiable ones while preserving injectivity of the homomorphism 
$G_\Ga\to Ham(S^2)$. In order to do so, we will need an estimate on the growth of partial derivatives of the functions $H_v$ at the unit circle. 

\begin{lem}\label{fast}
For every $v\in V(\Ga)$ and each $n=k+m$, the function $H=H_v$ satisfies
$$
\left|\frac{\partial^n}{\partial z^m \partial \bar{z}^k} H(z)\right|= O( r^{-(n-2)})
$$
where $r=1-|z|$. In particular, all $n$-th order derivatives of $H$ blow up at $S^1$ at most polynomially fast.  
\end{lem}
\proof The proof repeats the argument in Lemma \ref{lip}. Suppose that $z\in A_{v,\si}$. Set $w=\si^{-1}(z)\in A(v)$. 
Then $H(z)= \la_{v,\si}^{2} \hat{H}_{v,\si}\circ \si^{-1}(z)$, 
where
$$
\la_{v,\si}^{-1}\about \left|\frac{d}{dz} \si^{-1}(z)\right|.  
$$
The partial derivatives 
$$
\frac{\partial^k }{\partial \bar{z}^k} H(z)=   \la_{v,\si}^{2} \frac{\partial^k }{\partial \bar{z}^k}  \hat{H}_{v,\si}(w)
$$
are uniformly bounded (with respect to $\si$) for each $k$. On the other hand, since all derivatives of orders $\le k$ of the functions 
$\hat{H}_{v,\si}$ are uniformly bounded in $\si$, we have 
$$
\left|\frac{\partial^m}{ \partial z^m} H(z)\right|\about  \left| \la_{v,\si}^{2}\right|\cdot  \left|\frac{\partial^m }{\partial z^m} \si(w) \right|^{-1}. 
$$
Then we observe that
$$
\left|\frac{\partial^m }{\partial z^m} \si(w)\right| \about \la_{v,\si}^m\about r^m. 
$$
Lemma follows. \qed

\medskip
It is clear however that the above calculations cannot get bet\-ter than $C^{1,1}$-smooth\-ness for the function $H_v$.  In order to embed
$G_\Ga$ in $Ham(S^2)$ which consists of smooth Hamiltonian diffeomorphisms, we will use an approximation argument.

\section{Approximation}\label{approximation}

{\bf The mollifiers.} We define a family of $C^\infty$ functions $\eta_\eps(z), \eps>0$ (the mollifiers) on $S^2$ so that:

\begin{itemize}

\item For every $\eps>0, z\notin D$, $\eta_\eps(z)=0$. 
Moreover, $\eta_\eps(z)$ and its derivatives of all orders vanish exponentially fast on $S^1=\partial D$. 

\item $\max_{z\in D} \eta_\eps(z)=1= \eta_\eps(0)$. 

\item For every fixed $z\in D$ the functions 
$$
\eps\mapsto \eta_\eps(z), \quad \eps\mapsto d\eta_\eps(z)$$
are real-analytic. 

\item $$
\lim_{\eps\to 0} \eta_\eps(z)=1
$$
in $C^\infty$ topology uniformly on compacts in $D$. 

\end{itemize}

Explicitly, one can take 
$$
\eta_\eps(z)= \varphi_\eps(|z|), \quad |z|<1, 
$$
$$
\eta_\eps(z)=0, \quad |z|\ge 1. 
$$
where $\varphi_\eps(x)$ is the composition of 
$$
 \exp(-\eps y^2) 
$$
and
$$
y=\tan\left(\frac{\pi x}{2} \right). 
$$

Now, set $H_{v}^{(\eps)}:= \eta_\eps H_v$ where the functions $H_v$ on $S^2$ are defined by (\ref{H_v}). 
Since derivatives of all orders of $H_v$ blow up on $S^1$ at most polynomially fast 
(Lemma \ref{fast}), it follows that the functions $H_{v}^{(\eps)}$ are $C^\infty$ on $S^2$. Clearly, for $[uv]\notin E(\Ga)$, the supports of $H_{u}^{(\eps)}, H_{v}^{(\eps)}$ in $D$ are disjoint. Therefore, the corresponding Hamiltonian maps $f_{v,\eps}$
(with respect to the form $\om_0$)  commute. Moreover, the functions $H_{v}^{(\eps)}$ (and their derivatives) 
depend analytically on $\eps$ and converge to $H_v$ in $C^\infty$-topology uniformly on compacts in the open disk $D$. Therefore, the corresponding time-$N$ Hamiltonian maps $\rho_\eps(g_v):=f_{v,\eps}^N$   
converge to $\rho_0(g_v)$ as well (uniformly on compacts in $D$), where $\rho_0$ is defined by the formula (\ref{rho_0}). 
Moreover, for each $v$ and $z$, the function 
$$
\eps \mapsto f_{v,\eps}^N(z) 
$$
is real-analytic, for $\eps>0$.

We therefore obtain a family of representations $\rho_\eps: G_\Ga\to Ham(S^2)$ which send the generators $g_v$ to $\rho_\eps(g_v)$ as above.

\begin{lem}
For all but countably many $\eps$, the representations $\rho_\eps$ are faithful.  
\end{lem}
\proof For a fixed $g\in G_\Ga \setminus \{1\}$ the set $E_g$ of $\eps>0$ for which $g\in Ker(\rho_\eps)$ is either countable or the entire 
$\R_+$ (since $\rho_\eps(g)$ depends real-analytically on $\eps$). If all the sets $E_g$ are countable, we are done. Otherwise, 
there exists $g\in G_\Ga \setminus \{1\}$ which maps trivially by all $\rho_\eps$. Then
the limit
$$
\rho_0(g)= \lim_{\eps\to 0} \rho_\eps(g)
$$
is also the identity on $D$. However, this contradicts faithfulness of $\rho_0$. \qed 

This concludes the proof of Theorem \ref{main}. \qed 

\medskip
{\bf Higher-dimensional symplectic manifolds.}

{\em Proof of Corollary \ref{c:main}.} 
Let $2n$ be the dimension of $M$. Consider a polydisk $D^n\subset M$, where $D\subset \C$ is the unit disk, embedded in $M$ so that restriction of 
the symplectic structure $\om$ on $D^n$  splits as the sum
$$
c\cdot \om_0 \oplus ... \oplus \om_0
$$
where $\om_0$ is the Euclidean area form on each factor and $c$ is a sufficiently small positive constant. 
Take a faithful representation $\rho_\eps: G_\Ga\to Ham(D, \om_0)\subset Ham(S^2,\om_0)$ constructed in the proof of Theorem \ref{main}. Then the group $\rho_\eps(G_\Ga)$ fixes the boundary of $D$ pointwise. 
The images of the generators $\rho_\eps(g_v)$ are time-$N$ maps of functions $k_v:=H_{v}^{(\eps)}$ supported in $D$ (where $N\ge 2$). 
Then we define the function $h_v: D^n= D_1\times ...\times D_n \to \R$ by 
$$
h_v(z_1,...,z_n)= k_v(z_1) \eta(z_2)...\eta(z_n), 
$$
where $\eta(z):=\eta_1(z)$, see the definition of the mollifier $\eta_t$ in \S \ref{approximation}.  

\begin{lem}
The function $h_v$ vanishes on the boundary of  $D^n$ with all its derivatives. 
\end{lem}
\proof Let $p$ be a boundary point of $D^n$. If $p\in \partial D_1 \times D_2\times ... \times D_n$, the assertion follows 
from the fact that the function $k_v(z)$ vanishes with all its derivatives on the boundary circle of $D=D_1$. If 
$p\in D\times \partial (D_2\times ... \times D_n)$, then vanishing follows from vanishing of $\eta$ with all its derivatives at 
the boundary of $D$. \qed 

\medskip
We, thus, extend $h_v$ by zero to the rest of the manifold $M$ and retain the notation $h_v$ for the extension. Note that the  supports 
of $h_v, h_w$ in $D^n$ are disjoint provided that $[v,w]\notin E(\Ga)$. We next observe that at every point $z=(z_1,0,...,0)\in D_1\times 0 \times ... \times 0\subset D^n$, the differential of $h_v$ equals
$$
dh_v(z_1,0,...,0)= d k_v(z_1)
$$ 
since $d\eta(0)=0$ and $\eta(0)=1$. Therefore, the time-$N$ map $\rho(g_v)$ of $h_v$ is supported in the polydisk $D^n$ and satisfies 
$$
\rho(g_v): (z_1,0,...0) \mapsto (\rho_\eps(g_v)(z_1), 0,...0).  
$$
Hence, $\rho(g_v)$ preserves the disk $D_1\times  0 \times ... \times 0\subset D^n$ and acts on this disk as $\rho_\eps(g_v)$. 
It is then clear that the map $g_v\mapsto \rho(g_v)$ determines a monomorphism $G_\Ga\to Ham(M,\om)$, since 
$\rho_\eps: G_\Ga\to Ham(D)$ was injective. \qed

Address:

Michael Kapovich: Department of Mathematics, University of California, Davis, CA 95616,
USA. (kapovich@math.ucdavis.edu)

\end{document}